\documentclass[12pt]{article}
\usepackage{latexsym}
\usepackage{amsmath,amssymb}
\usepackage{mathrsfs}
\usepackage{amsthm}
\usepackage{amsfonts}
\numberwithin{equation}{section}
\newtheorem{definition}{Definition}[section]
\newtheorem{theorem}{Theorem}[section]
\newtheorem{lemma}{Lemma}[section]

\newtheorem{proposition}{Proposition}[section]
\newtheorem{hypothesis}{Hypothesis}

\newtheorem{corollary}{Corollary}[section]

\begin{document}
\author{Hanzhong Wu \\School of Mathematical Sciences and Key Lab of Mathematics\\
for Nonlinear Sciences, Fudan University, Shanghai 200433,
China\\Email: hzwu@fudan.edu.cn}
\date{}
\title{On Certain Hypotheses in Optimal Control Theory and the Relationship of the Maximum
Principle with the Dynamic Programming Method\\ Proposed by L. I.
Rozonoer } \maketitle \noindent{{\bf Abstract.} In this paper we
will study three hypotheses proposed by L. I. Rozonoer \cite{R} in
optimal control theory in order to derive conditions for the
existence of an optimal control under all initial conditions, and
the relationships between Pontryagin maximum principle and the
dynamic programming method.
}\\

\section{Introduction}

Let us introduce the following optimal control problem considered in
\cite{R}:

\vskip 2mm

\noindent{\bf OCP} {\it To minimize the Lagrange cost functional
\begin{eqnarray}
\int_{t_0}^TF(x,u,t)\,dt
\end{eqnarray}
subject to the controlled system
\begin{eqnarray}
\dot{x}=f(x,u,t),
\end{eqnarray}
with $u(t)\in U$ and the initial state condition
\begin{eqnarray}
x(t_0)=x^0,
\end{eqnarray}
where $t_0$ and $T$ with $t_0<T$ are prescribed real numbers.}

In {\bf OCP}, $U\subseteq \mathbb{R}^m$ is the control domain while
the set of admissible controls under consideration is the set of all
Lebesque measurable selection $u(t)\in U$ ( see also (2.12) or
Remark 4.2); $x^0=(x_1^0,\cdots,x_n^0)^T\in\mathbb{R}^n$ is the
initial state, $x=(x_1,\cdots,x_n)^T\in\mathbb{R}^n$ is the state
variable, and $f=(f_1,\cdots,f_n)^T$ is n-dimensional vector-valued
function, where and throughout this paper the superscript $^T$
denotes the transpose of a vector or matrix. Other technical
assumptions on $f$ and $F$ will be given in the following sections.

The control Hamiltonian for {\bf OCP} is
\begin{eqnarray}
\mathscr{H}(x,p,u,t):=\sum_{i=1}^np_if_i(x,u,t)-F(x,u,t),
\end{eqnarray}
where $p=(p_1,\cdots,p_n)^T\in \mathbb{R}^n$ is the costate
variable.

For any given initial data $(x^0,\tau)$ with $\tau\in [t_0,T)$ and
$x^0\in \mathbb{R}^n$, we introduce the control Hamiltonian system
\begin{eqnarray}
\left\{
\begin{array}{l}
\dot{x}=\frac{\partial \mathscr{H}(x,p,u,t)}{\partial p},\\
\dot{p}=-\frac{\partial \mathscr{H}(x,p,u,t)}{\partial x},
\end{array}\right.
\end{eqnarray}
with the two-point boundary value conditions
\begin{eqnarray}
x(\tau)=x^0,\qquad p(T)=0.
\end{eqnarray}

\noindent{\bf Remark 1.1.}\quad In order to distinguish the function
(1.4) and the system (1.5) with the Hamiltonian (4.10) and the
canonical Hamiltonian system (4.21), which will be considered in
Section 4 and very related to these analogues, we prefer to calling
(1.4) (and (1.5)) the control Hamiltonian (and the control
Hamiltonian system) instead of the Hamiltonian (and the Hamiltonian
system).

\begin{definition} A control $u^*(\cdot):[\tau,T]\mapsto U$ is said to satisfy
the Pontryagin maximum condition on the interval $[\tau, T]$ under
the initial data $(x^0,\tau)$, provided that the unique solution
$(x^*(\cdot),p^*(\cdot))$ of the control Hamiltonian system
(1.5)-(1.6) corresponding to this control $u^*(\cdot)$ satisfy
\begin{eqnarray}
\mathscr{H}(x^*(t),p^*(t),u^*(t),t)\ge
\mathscr{H}(x^*(t),p^*(t),u,t),\qquad\forall u\in U,\quad\text{a.e.
}t\in [\tau,T].
\end{eqnarray}
\end{definition}

Related to {\bf OCP}, the Hamilton-Jacobi-Bellman equation (or
called the Bellman equation in \cite{R}) is
\begin{eqnarray}
-v_{\tau}+\sup_{u\in U}\mathscr{H}(x,-v_{x},u,\tau)=0,\qquad
(x,\tau)\in \mathbb{R}^n\times (t_0,T),
\end{eqnarray}
with the boundary condition
\begin{eqnarray}
v(x,T)=0,
\end{eqnarray}
where and throughout this paper, the partial derivative of a given
function $\varphi$ with respect to $\tau\in [t_0,T]$ or
$x\in\mathbb{R}^n$ will be denoted by $\varphi_{\tau}$ or
$\varphi_x$, respectively.

\vskip 2mm

In order to adapt for the optimal control theory, L. I. Rozonoer
\cite{R} first give the following concept of weak solution to the
Hamilton-Jacobi-Bellman equation:

\begin{definition} A continuous function $V:\mathbb{R}^n\times [t_0,T]\mapsto \mathbb{R}$ is said to be an extended
solution of the Hamilton-Jacobi-Bellman equation (1.8)-(1.9)
provided that, for any $x\in \mathbb{R}^n$ and $\tau\in [t_0,T]$
there exists a $p^*$ with $-p^*\in\partial_+V_x(x,\tau)$ and a
$u^*\in U$ such that
\begin{eqnarray}
\mathscr{H}(x,p^*,u^*,\tau)\ge
\mathscr{H}(x,p^*,u,\tau),\qquad\forall u\in U,
\end{eqnarray}
and
\begin{eqnarray}
\mathscr{H}(x,p^*,u^*,\tau)\in\partial_+V_{\tau}(x,\tau),
\end{eqnarray}
along with the boundary condition (1.9).
\end{definition}
The notation $\partial_+V_x(x,\tau)$ and
$\partial_+V_{\tau}(x,\tau)$ in {\bf Definition 1.2} means the
superdifferential of the function $V(\cdot,\cdot)$ at the point
$(x,\tau)$ with respect to $x$ and $\tau$, respectively. We will
recall the concepts of superdifferential and subdifferential in
Section 2.

\vskip 2mm

In order to derive conditions for the existence of an optimal
control under all initial conditions, and thereby the relationships
between Pontryagin maximum principle and the dynamic programming
method, L. I. Rozonoer \cite{R} proposed three hypotheses on {\bf
OCP}.

\begin{hypothesis} The existence of an extended solution to the Hamilton-Jacobi-Bellman
equation is necessary and sufficient for the existence of an optimal
control under all initial data $(x^0,\tau)$.
\end{hypothesis}

\begin{hypothesis} The extended solution of the Hamilton-Jacobi-Bellman
equation exists if and only if for every initial data $(x^0,\tau)$,
there exists a unique control satisfying the Pontryagin maximum
condition.
\end{hypothesis}

\begin{hypothesis} If for every initial data $(x^0,\tau)$, there
exists a unique control satisfying the Pontryagin maximum condition,
then this control is optimal.
\end{hypothesis}

Just as emphasized in \cite{R,R1} that the concept of solution to
the Hamilton-Jacobi-Bellmen equation need to be generalized to
ensure that certain general hypotheses could be given on the
condition for the existence of an optimal control under all initial
data, and as a result the relationships between Pontryagin maximum
principle and the dynamic programming method. In this approach,
\cite{K1,K2} generalized the concept of solution of the
Hamilton-Jacobi-Bellmen equation to help demonstrating the necessary
and sufficient conditions for minimization of a functional not only
for nonsmooth cases, but also for the case where there is even no
optimal control, and provide a possibility for investigating control
design. On the other hand, many works such as \cite{XYZ, YZ} and the
references cited within devoted to the relationships between
Pontryagin maximum principle and the dynamic programming method
directly or in the framework of viscosity solution theory of the
Hamilton-Jacobi-Bellmen equation. There are rich references related
to this approach (see \cite{R,YZ}).

In this paper, we will only focus on the problem {\bf OCP} in order
to study these above hypotheses.

The rest of paper is organized as follows: In Section 2, we will
study {\bf Hypothesis 1}. First, it will be considered the concept
of the extended solution defined by L. I. Rozonoer. Under some mild
technical assumptions, the Bellman function is just right an
extended solution to the Hamilton-Jacobi-Bellmen equation. Second,
one example will be given to show that the Bellman function is not
an extended solution but a viscosity solution to the
Hamilton-Jacobi-Bellmen equation, which indicates the application
range of the extended solution in some sense. Finally, one
counterexample will be given to verify that {\bf OCP} may have no
optimal controls under some initial data $(x^0,\tau)$ even that the
Bellman function is a classical ($C^2$ smooth) solution to the
Hamilton-Jacobi-Bellmen equation, which is verified to be an
extended solution as well. In Section 3, we will study {\bf
Hypothesis 2}. Two counterexamples will be given to show that there
are many optimal controls under every initial data $(x^0,\tau)$ even
that the Bellman function is a classical ($C^2$ smooth) solution to
the Hamilton-Jacobi-Bellmen equation, which is verified to be an
extended solution as well. On the other hand, all optimal controls
satisfy the Pontryagin maximum condition. In section 4, we will
study {\bf Hypothesis 3}. First, it will be given the necessary and
sufficient condition for the differentiability of the Hamiltonian.
Then, main results will be established that {\bf Hypothesis 3} holds
true under some technical assumptions of regularity on the data,
through the existing relationships between the Hamiltonian system
and the Hamilton-Jacobi equation. In Section 5, the conclusions will
be given.

\section{On the extended solution and Hypotheses 1}

\subsection{On the concept of the extended solution}

First, we recall some related concepts and results in the theory of
the Hamilton-Jacobi-Bellman equations.

Let $n$ be a positive integer. We denote by the operation
$\langle\cdot,\cdot\rangle$ and $\|\cdot\|$ the inner product and
norm in $\mathbb{R}^n$.

The following definition is combined from \cite{BC}, \cite{CS} and
\cite{YZ}, etc.

\begin{definition} For a continuous function $\varphi:\mathbb{R}^n\times [t_0,T]\mapsto
\mathbb{R}$, define that
\begin{eqnarray}
\begin{array}{l}
\partial_+\varphi_{x,\tau}(x,\tau):=\{(p,q)\in\mathbb{R}^{n+1}|\limsup_{(y,t)\rightarrow (x,\tau)}\frac{\varphi(y,t)-\varphi(x,\tau)-\langle p,y-x\rangle-q(t-\tau)}{|t-\tau|+\|y-x\|}\le
0\},\\
\partial_-\varphi_{x,\tau}(x,\tau):=\{(p,q)\in\mathbb{R}^{n+1}|\liminf_{(y,t)\rightarrow (x,\tau)}\frac{\varphi(y,t)-\varphi(x,\tau)-\langle p,y-x\rangle-q(t-\tau)}{|t-\tau|+\|y-x\|}\ge
0\},\\
\partial_+\varphi_{x}(x,\tau):=\{p\in\mathbb{R}^n|\limsup_{y\rightarrow x}\frac{\varphi(y,\tau)-\varphi(x,\tau)-\langle p,y-x\rangle}{\|y-x\|}\le
0\},\\
\partial_-\varphi_{x}(x,\tau):=\{p\in\mathbb{R}^n|\liminf_{y\rightarrow x}\frac{\varphi(y,\tau)-\varphi(x,\tau)-\langle p,y-x\rangle}{\|y-x\|}\ge
0\},\\
\partial_+\varphi_{\tau}(x,\tau):=\{q\in\mathbb{R}|\limsup_{t\rightarrow \tau}\frac{\varphi(x,t)-\varphi(x,\tau)-q(t-\tau)}{|t-\tau|}\le
0\},
\end{array}
\end{eqnarray}
for a given $(x,\tau)\in\mathbb{R}^n\times (t_0,T)$.

$\partial_+\varphi_{x,\tau}(x,\tau)$ and
$\partial_-\varphi_{x,\tau}(x,\tau)$ are called the
superdifferential and subdifferential of $\varphi$ at $(x,\tau)$,
respectively; $\partial_+\varphi_x(x,\tau)$ and
$\partial_-\varphi_x(x,\tau)$ are called the partial
superdifferential and subdifferential of $\varphi$ at $(x,\tau)$
with respect to $x$, respectively;
$\partial_+\varphi_{\tau}(x,\tau)$ is called the partial
superdifferential of $\varphi$ at $(x,\tau)$ with respect to
$\tau$.
\end{definition}

\noindent{\bf Remark 2.1.}\quad We can define the right
superdifferential $\partial_+\varphi_{x,\tau+}(x,\tau)$, right
subdifferential $\partial_-\varphi_{x,\tau+}(x,\tau)$ and partial
right superdifferential $\partial_+\varphi_{\tau+}(x,\tau)$ with
respect to $\tau$ at $(x,\tau)\in\mathbb{R}^n\times [t_0,T)$ by
restricting $t\downarrow\tau$ in (2.1). Analogously, define the left
superdifferential $\partial_+\varphi_{x,\tau-}(x,\tau)$, left
subdifferential $\partial_-\varphi_{x,\tau-}(x,\tau)$ and partial
left superdifferential $\partial_+\varphi_{\tau-}(x,\tau)$ at
$(x,\tau)\in\mathbb{R}^n\times (t_0,T]$ by restricting
$t\uparrow\tau$ in (2.1).

\vskip 4mm

From the above definition, it can be easily deduced that
\begin{lemma} Let $\varphi:\mathbb{R}^n\times [t_0,T]\mapsto \mathbb{R}$ be
continuous. It holds that
\begin{description}
\item[(a)]
if $(p,q)\in \partial_+\varphi_{x,\tau}(x,\tau)$, then $p\in
\partial_+\varphi_x(x,\tau)$ and $q\in
\partial_+\varphi_{\tau}(x,\tau)$;
\item[(b)]
if $(p,q)\in \partial_+\varphi_{x,\tau+}(x,t_0)$, then $p\in
\partial_+\varphi_x(x,t_0)$ and $q\in
\partial_+\varphi_{\tau+}(x,t_0)$;
\item[(c)]
if $(p,q)\in \partial_+\varphi_{x,\tau-}(x,T)$, then $p\in
\partial_+\varphi_x(x,T)$ and $q\in
\partial_+\varphi_{\tau-}(x,T)$.
\end{description}
\end{lemma}

The definition of the viscosity solution to the first order PDEs is
first given by Crandall and Lions \cite{CP}. We can also refer to
\cite{BC}, \cite{CS}, \cite{F} and \cite{YZ}, etc., for the
following definition.

\begin{definition} A continuous function $\varphi:\mathbb{R}^n\times [t_0,T]\mapsto \mathbb{R}$ is called a viscosity subsolution of the
Hamilton-Jacobi-Bellman equation (1.8) provided that, for any
$(x,\tau)\in\mathbb{R}^n\times (t_0,T)$,
\begin{eqnarray}
-q+\sup_{u\in U}\mathscr{H}(x,-p,u,\tau)\le 0,\qquad\forall
(p,q)\in\partial_+\varphi_{x,\tau}(x,\tau);
\end{eqnarray}
A continuous function $\varphi:\mathbb{R}^n\times [t_0,T]\mapsto
\mathbb{R}$ is called a viscosity supersolution of the
Hamilton-Jacobi-Bellman equation (1.8) provided that, for any
$(x,\tau)\in\mathbb{R}^n\times (t_0,T)$,
\begin{eqnarray}
-q+\sup_{u\in U}\mathscr{H}(x,-p,u,\tau)\ge 0,\qquad\forall
(p,q)\in\partial_-\varphi_{x,\tau}(x,\tau).
\end{eqnarray}
Finally, $\varphi$ is called a viscosity solution of the
Hamilton-Jacobi-Bellman equation (1.8) if it is simultaneously a
viscosity sub- and supersolution. In addition, if $\varphi$
satisfies the boundary condition (1.9), then $\varphi$ is called a
viscosity solution of the Hamilton-Jacobi-Bellman equation
(1.8)-(1.9).
\end{definition}

Denote $\mathbb{R}_+=[0,+\infty)$. For any $R>0$, we denote by
$B_R(\mathbb{R}^{n+1})$ (or $B_R(\mathbb{R}^n)$) the open ball in
$\mathbb{R}^{n+1}$ (or $\mathbb{R}^n$) with a radius $R$ centered at
$0$.
\begin{definition} Consider a convex subset $K\subseteq \mathbb{R}^{n+1}$ (or $\mathbb{R}^n$). A function $\varphi:
K\mapsto \mathbb{R}$ is called semiconcave if there exists a
function $\omega: \mathbb{R}_+\times \mathbb{R}_+\mapsto
\mathbb{R}_+$ satisfies that
\begin{eqnarray}\left\{\begin{array}{ll}\omega(r,d)\le\omega(R,D),&\quad\forall r\le R, d\le D,\\
\lim_{D\rightarrow 0+}\omega(R,D)=0,&\quad\forall
R>0,\end{array}\right.
\end{eqnarray}
such that, for every $R>0$, $\lambda\in [0,1]$ and any $\xi,\eta\in
K\cap B_R(\mathbb{R}^{n+1})$ (or $K\cap B_R(\mathbb{R}^n)$),
\begin{eqnarray}
\lambda\varphi(\xi)+(1-\lambda)\varphi(\eta)-\varphi[\lambda
\xi+(1-\lambda)\eta]\le\lambda(1-\lambda)\|\xi-\eta\|\omega(R,\|\xi-\eta\|).
\end{eqnarray}
We call the above function $\omega$ a modulus of a semiconcavity of
$\varphi$.
\end{definition}

\noindent{\bf Remark 2.2.} \quad Obviously, if
$\varphi:\mathbb{R}^n\times [t_0,T]\mapsto \mathbb{R}$ is
semiconcave, then both
$\varphi(\cdot,t_0):\mathbb{R}^n\mapsto\mathbb{R}$ and
$\varphi(\cdot,T):\mathbb{R}^n\mapsto\mathbb{R}$ are semiconcave.

\vskip 3mm

By Rademacher's theorem (\cite{E}, Ch.5, p.281), if
$\varphi:\mathbb{R}^n\times [t_0,T]\mapsto \mathbb{R}$ is locally
Lipschitz continuous, then $\varphi$ is differentiable almost
everywhere in $\mathbb{R}^n\times [t_0,T]$. Meanwhile,
\begin{eqnarray*}
D^*\varphi(x,\tau):=&&\{\lim_{i\rightarrow+\infty}\varphi_{x,\tau}(x_i,\tau_i)\hskip
1mm |\hskip 1mm \mathcal{D}(\varphi)\ni (x_i,\tau_i)\rightarrow
(x,\tau)\\&&\hskip 2mm \text{ such that
}\lim_{i\rightarrow+\infty}\varphi_{x,\tau}(x_i,\tau_i) \text{
exists}\},
\end{eqnarray*}
is nonempty at all $(x,\tau)\in \mathbb{R}^n\times [t_0,T]$, where
$$
\mathcal{D}(\varphi):=\{(x,\tau)\in \mathbb{R}^n\times (t_0,T)\hskip
1mm |\hskip 1mm \varphi\text{ is differentiable at }(x,\tau)\}.
$$

\begin{lemma} Let $\varphi:\mathbb{R}^n\times [t_0,T]\mapsto \mathbb{R}$ be locally
Lipschitz continuous. If $\varphi$ is semiconcave, then it holds
that
\begin{description}
\item[(a)] for all $(x,\tau)\in \mathbb{R}^n\times (t_0,T)$,
$$
\partial_+\varphi_{x,\tau}(x,\tau)={\rm co} D^*\varphi(x,\tau)\ne\emptyset,
$$
where the operation `` {\rm co}" denotes the convex hull;
\item[(b)]
$$\partial_+\varphi_{x,\tau+}(x,t_0)\supseteq{\rm co}
D^*\varphi(x,t_0)\ne\emptyset,$$ and
$$\partial_+\varphi_{x,\tau-}(x,T)\supseteq {\rm
co}D^*\varphi(x,T)\ne\emptyset. $$
\end{description}
\end{lemma}

{\bf Proof}\quad  Part (a) is one part of Theorem 3.3.6 in \cite{CS}
or in \cite{BC}. Part (b) follows easily from the locally Lipschitz
continuity and the semiconcavity of $\varphi$ on $\mathbb{R}^n\times
[t_0,T]$, similar to the proofs of Proposition 3.3.1 and 3.3.4 in
\cite{CS}.\hfill$\Box$\vspace{3mm}

\noindent{\bf Remark 2.3.}\quad The semiconcavity of
$\varphi:\mathbb{R}^n\times [t_0,T]\mapsto \mathbb{R}$ implies the
locally Lipschitz continuity only on $\mathbb{R}^n\times (t_0,T)$.
(see \cite{CS,BC})

\vskip 3mm

In this section, we will need some technical assumptions on
$f:\mathbb{R}^n\times U\times [t_0,T]\mapsto \mathbb{R}^n$ and
$F:\mathbb{R}^n\times U\times [t_0,T]\mapsto \mathbb{R}$ as follows:

\begin{description}
\item[(H1)] Both $f$ and $F$
are continuous, and there exists a constant $M>0$ such that
\begin{eqnarray}
\|f(0,u,t)\|\le M,\qquad \forall (u,t)\in U\times [t_0,T],
\end{eqnarray}
and
\begin{eqnarray}
|F(0,u,t)|\le M,\qquad \forall (u,t)\in U\times [t_0,T].
\end{eqnarray}
\item[(H2)] Both $f(\cdot,u,\cdot)$ and $F(\cdot,u,\cdot)$ are locally Lipschitz
continuous on $\mathbb{R}^n\times [t_0,T]$ uniformly in $u\in U$,
i.e., there exists a nondecreasing function
$L:\mathbb{R}_+\mapsto\mathbb{R}_+$ such that, for any $(x,t),
(y,s)\in B_R(\mathbb{R}^n)\times [t_0,T]$, and any $u\in U$,
\begin{eqnarray}
\|f(x,u,t)-f(y,u,s)\|\le L(R)(|t-s|+\|x-y\|),
\end{eqnarray}
and
\begin{eqnarray}
|F(x,u,t)-F(y,u,s)|\le L(R)(|t-s|+\|x-y\|).
\end{eqnarray}
\item[(H3)] There exists a modulus $\omega: \mathbb{R}_+\times \mathbb{R}_+\mapsto
\mathbb{R}_+$ satisfying (2.4) such that, for any $\lambda\in [0,1]$
and $(x,t), (y,s)\in B_R(\mathbb{R}^n)\times [t_0,T]$, and any $u\in
U$,
\begin{eqnarray}\begin{array}{l}
\|\lambda f(x,u,t)+(1-\lambda)f(y,u,s)-f[\lambda
x+(1-\lambda)y,u,\lambda
t+(1-\lambda)s]\|\\\le\lambda(1-\lambda)(|t-s|+\|x-y\|)\omega(R,|t-s|+\|x-y\|)
.\end{array}
\end{eqnarray}
\item[(H4)] $F(\cdot,u,\cdot)$ is semiconcave on $\mathbb{R}^n\times [t_0,T]$, i.e., there exists a modulus $\omega: \mathbb{R}_+\times \mathbb{R}_+\mapsto
\mathbb{R}_+$ satisfying (2.4) such that, for any $\lambda\in [0,1]$
and $(x,t), (y,s)\in B_R(\mathbb{R}^n)\times [t_0,T]$, and any $u\in
U$,
\begin{eqnarray}\begin{array}{l}
\lambda F(x,u,t)+(1-\lambda)F(y,u,s)-F[\lambda
x+(1-\lambda)y,u,\lambda
t+(1-\lambda)s]\\\le\lambda(1-\lambda)(|t-s|+\|x-y\|)\omega(R,|t-s|+\|x-y\|)
.\end{array}
\end{eqnarray}
\end{description}
\noindent{\bf Remark 2.4.}\quad {\bf (H3)} holds true in particular
when $f$ is continuously differentiable with respect to $(x,t)$
uniformly in $u$. More precisely, if we assume that there exists a
modulus $\omega: \mathbb{R}_+\times \mathbb{R}_+\mapsto
\mathbb{R}_+$ satisfying (2.4) such that, for any $(x,t), (y,s)\in
B_R(\mathbb{R}^n)\times [t_0,T]$, and any $u\in U$,
$$
\|f_t(x,u,t)-f_t(y,u,s)\|+\|f_x(x,u,t)-f_x(y,u,s)\|\le\omega(R,|t-s|+\|x-y\|).
$$
Conversely, under the assumption {\bf (H2)}, it follows from
Proposition 1.1.13 in \cite{F} that, {\bf (H3)} implies that $f$ is
continuously differentiable with respect to $(x,t)$.

\vskip 3mm

Throughout this paper, we define the set of all admissible controls
under any given initial time $\tau\in [t_0,T]$, as follows:
\begin{eqnarray}
{\cal U}(\tau):=\{u(\cdot):[\tau,T]\mapsto U|\quad u(\cdot)\text{ is
Lebesgue measurable}\},
\end{eqnarray}
and denote by ${\cal U}$ simply for ${\cal U}(t_0)$. For any given
initial data $(x^0,\tau)\in\mathbb{R}^n\times [t_0,T]$,
$x(\cdot;x^0,\tau,u(\cdot))$ is the unique solution of the control
system (1.2) under the initial state condition $x(\tau)=x^0$ and the
control $u(\cdot)\in{\cal U}(\tau)$, which we will only denote by
$x(\cdot)$ for short if without confusion. The value function (or
called the Bellman function in \cite{R})
$V(\cdot,\cdot):\mathbb{R}^n\times [t_0,T]\mapsto\mathbb{R}$ is
defined by
\begin{eqnarray}
V(x^0,\tau):=\inf_{u(\cdot)\in{\cal
U}(\tau)}\int_{\tau}^TF(x(t),u(t),t)\,dt,
\end{eqnarray}
where $x(\cdot)$ is the solution of (1.2) under the initial state
condition $x(\tau)=x^0$ and the control $u(\cdot)\in{\cal U}(\tau)$.
\vskip 3mm

Similar to Theorem 4.1 in \cite{Ca} (or \cite{F, CS} etc.), it
follows that
\begin{proposition} Assume that (H1)-(H2) are satisfied. Then the value function $V(\cdot,\cdot)$
is locally Lipschitz continuous on $\mathbb{R}^n\times [t_0,T]$.
\end{proposition}

Similar to Theorem 3.2 in \cite{Ca} (or \cite{F, CS} etc.), it
follows that
\begin{proposition} Assume that (H1)-(H4) are satisfied. Then the value function $V(\cdot,\cdot)$
is semiconcave on $\mathbb{R}^n\times [t_0,T]$.
\end{proposition}

\noindent{\bf Remark 2.5.} \quad The technical assumptions (H1)-(H2)
 in Proposition 2.1 and (H1)-(H4)
 in Proposition 2.2 can be weaken in some approaches, for example, the consideration of the cases with unbounded control variables in \cite{Ca},
 etc.

\begin{theorem} Assume that (H1)-(H4) are satisfied. If the control domain $U$ is compact, then the value function
$V(\cdot,\cdot):\mathbb{R}^n\times [t_0,T]\mapsto\mathbb{R}$ is an extended
solution of the Hamilton-Jacobi-Bellman equation (1.8)-(1.9).
\end{theorem}
{\bf Proof}\quad It is well known that, $V$ is the unique viscosity
solution to the Hamilton-Jacobi-Bellman equation
(1.8)-(1.9).

Proposition 2.1 and 2.2 yields that the value function $V$ is
locally Lipschitz continuous and semiconcave on $\mathbb{R}^n\times
[t_0,T]$.

By Rademacher's theorem (\cite{E}, Ch.5, p.281), the locally
Lipschitz continuity of the value function $V$ implies that $V$ is
differentiable almost everywhere in $\mathbb{R}^n\times (t_0,T)$.
Hence, for any given $(x,\tau)\in\mathbb{R}^n\times [t_0,T]$, there
exists a sequence of $\{(x_i,\tau_i)\}\subset\mathbb{R}^n\times
(t_0,T)$ such that
$$
(x_i,\tau_i)\rightarrow (x,\tau),\qquad \text{as}\quad
i\rightarrow+\infty,
$$
and $V$ is differentiable at all $(x_i,\tau_i)$. According to
Proposition 1.9 in \cite{BC} (or Theorem 1 in \cite{E}, Ch.10,
p.545), the value function $V$ satisfies the Hamilton-Jacobi-Bellman
equation (1.8) at $(x_i,\tau_i)$ in the classical sense, i.e.,
\begin{eqnarray}
-V_{\tau}(x_i,\tau_i)+\sup_{u\in
U}\mathscr{H}(x_i,-V_{x}(x_i,\tau_i),u,\tau_i)=0.
\end{eqnarray}
Meanwhile, due to the compactness of $U$, there exists a sequence of
$\{u_i\}\subset U$ such that
\begin{eqnarray}
-V_{\tau}(x_i,\tau_i)+\mathscr{H}(x_i,-V_{x}(x_i,\tau_i),u_i,\tau_i)=0,
\end{eqnarray}
and there exists a subsequence of $\{u_i\}$ (still denoted by
themselves without loss of generality) such that
$$u_i\rightarrow u^*\in U,\qquad \text{as}\quad
i\rightarrow+\infty.$$
Denote that
$$p_i:=-V_{x}(x_i,\tau_i),\qquad q_i:=V_{\tau}(x_i,\tau_i).$$
It follows from Lemma 2.2 that, there exists a subsequence of
$\{(-p_i,q_i)\}$ (still denoted by themselves without loss of
generality) such that
\begin{eqnarray}
V_{x,\tau}(x_i,\tau_i)\equiv(-p_i,q_i)\rightarrow
(-p^*,q^*)\in\left\{\begin{array}{l}\partial_+V_{x,\tau}(x,\tau),\\\partial_+V_{x,\tau+}(x,t_0),\\\partial_+V_{x,\tau-}(x,T),\end{array}\right.
\begin{array}{l}\text{if }\tau\in (t_0,T),\\\text{if }\tau=t_0,\\\text{if }\tau=T,\end{array}
\end{eqnarray}
as $i\rightarrow+\infty$, where $V_{x,\tau}(x_i,\tau_i)$ is the
derivative of $V(\cdot,\cdot)$ with respect to $(x,\tau)$ at
$(x_i,\tau_i)$.

By Lemma 2.1, combining (2.15) and (2.16) yields the conclusions.
\hfill$\Box$\vspace{3mm}

\subsection{An example for the concept of the extended solution}

\vskip 3mm

\noindent{{\bf Example.} Let the control domain be $U=[0,1]$. For
any given initial data $(x^0,\tau)$ with $\tau\in [0,1)$ and $x^0\in
\mathbb{R}$. Consider the following linear one-dimensional control
system
\begin{eqnarray}
\frac{dx}{dt}=x(t)+u(t), \qquad t\in [\tau,1],
\end{eqnarray}
with the initial condition
\begin{eqnarray}
x(\tau)=x^0\in \mathbb{R},
\end{eqnarray}
and let the associated cost functional be
\begin{eqnarray}
J(x^0,\tau;u(\cdot))=\int_{\tau}^1|x(t)|\,dt,
\end{eqnarray}
where the set of all admissible controls is
\begin{eqnarray}
{\cal U}(\tau):=\{u(\cdot):[\tau,1]\mapsto U|\quad u(\cdot)\text{ is
Lebesgue measurable}\}.
\end{eqnarray}

Obviously, this optimal control problem has a unique optimal control

\begin{eqnarray}
u^*(t)\equiv\left\{\begin{array}{ll}0,&\text{
 if } x^0\ge 0,\\1,&
\text{
 if } x^0<0, (1+x^0)e^{1-\tau}\le 1,\end{array}\right.
\end{eqnarray}
while
\begin{eqnarray}
u^*(t)\equiv\left\{\begin{array}{ll}1,&t\in
[0,\tau-\ln(1+x^0)),\\0,& t\in [\tau-\ln(1+x^0),1]
,\end{array}\right.
\end{eqnarray}
if $x^0<0$ and $(1+x^0)e^{1-\tau}> 1$.

The corresponding Hamilton-Jacobi-Bellman equation is
\begin{eqnarray}\left\{\begin{array}{l}
-\frac{\partial v}{\partial \tau}+\sup_{u\in U}\{(x+u)\frac{\partial
v}{\partial x}-|x|\}=0,\qquad (x,\tau)\in \mathbb{R}\times
(0,1),\\v(x,1)=0.\end{array}\right.
\end{eqnarray}

Obviously, the Hamilton-Jacobi-Bellman equation (2.23) has a unique
viscosity solution
\begin{eqnarray}
V(x,\tau)=\left\{\begin{array}{ll}x(e^{1-\tau}-1),&\text{
 if } x\ge 0,\\2+x-\tau-(1+x)e^{1-\tau},&
\text{
 if } x<0, (1+x)e^{1-\tau}\le 1,\\x-\ln(1+x),&
\text{
 if } x<0, (1+x)e^{1-\tau}> 1,\end{array}\right.
\end{eqnarray}
which is just the value function.

For any given $\tau_0\in (0,1)$, we have
\begin{eqnarray}
V(x,\tau_0)=\left\{\begin{array}{ll}x(e^{1-\tau_0}-1),&\text{
 if } x\ge 0,\\x-\ln(1+x),&
\text{
 if } e^{\tau_0-1}-1<x<0,\end{array}\right.
\end{eqnarray}
which implies that $\partial_+ V_x(0,\tau_0)=\emptyset$ and
$\partial_+ V_{x,\tau}(0,\tau_0)=\emptyset$.

Hence $V(\cdot,\cdot)$ defined by (2.24) is a viscosity solution but
not an extended solution of the Hamilton-Jacobi-Bellman equation
(2.23).

We notice that $V(\cdot,\cdot)$ is not a semiconcave function
according to Lemma 2.2.

\subsection{An example for Hypotheses 1}

We consider the following example of optimal control problem, which
is adapted from \cite{W} (Ch.3, p.246).

\vskip 3mm \noindent{{\bf Example.} Let the control domain be
$U=[-1,1]$, and the set of all control variables be
\begin{eqnarray}
{\cal U}:=\{u(\cdot):[0,1]\mapsto U|\quad u(\cdot)\text{ is Lebesgue
measurable}\}.
\end{eqnarray}
Consider the one-dimensional control system
\begin{eqnarray}
\frac{dx}{dt}=u(t), \qquad t\in [0,1],
\end{eqnarray}
with the initial state condition
\begin{eqnarray}
x(0)=x^0\in \mathbb{R},
\end{eqnarray}
and let the associated Lagrange type cost functional be
\begin{eqnarray}
J(x_0;u(\cdot))=\int_0^1[x^2(t)-u^2(t)]\,dt.
\end{eqnarray}

For this optimal control problem, the corresponding
Hamilton-Jacobi-Bellman equation is
\begin{eqnarray}\left\{\begin{array}{l}
-\frac{\partial v}{\partial \tau}+|\frac{\partial v}{\partial
x}|-x^2+1=0,\qquad (x,\tau)\in \mathbb{R}\times
(0,1),\\v(x,1)=0.\end{array}\right.
\end{eqnarray}

It is easy to verify that, the Hamilton-Jacobi-Bellman equation
(2.30) admits a $C^2$ solution
\begin{eqnarray}
V(x,\tau)=\left\{\begin{array}{ll}\frac{1}{3}[x^3-(x+\tau-1)^3]+\tau-1,&\text{
 if } x+\tau\ge 1,\\\frac{1}{3}x^3+\tau-1,&
\text{
 if } x+\tau<1,x\ge 0,\\-\frac{1}{3}x^3+\tau-1,&\text{
 if } -x+\tau<1,x<0,\\\frac{1}{3}[-x^3-(-x+\tau-1)^3]+\tau-1,&\text{
 if } -x+\tau\le
-1,\end{array}\right.
\end{eqnarray}
which is just the value function related to this optimal control
problem. Certainly, this solution is also an extended solution of
the HJB equation since the control Hamiltonian is
\begin{equation}\mathscr{H}(x,p,u,t)=pu-x^2+u^2,
\end{equation}
and
\begin{equation}\max_{u\in
[-1,1]}\mathscr{H}(x,p,u,t)=|p|-x^2+1,
\end{equation}
is attainable at $u=-1$ or $u=1$.

However, it can be proved similarly to \cite{W} (Ch. 3, p.247) that
there exists no optimal control under the initial data $(0,\tau)$
with $\tau\in [0,1)$. In these cases, $V(0,\tau)=-1+\tau$, which is
not attainable. In fact, there exists no optimal control under any
initial data $(x^0,\tau)$ with $\tau\in [0,1)$ and $|x^0|+\tau<1$.

\section{On Hypotheses 2}

\subsection{The first example}

Consider the special cases of {\bf OCP} with the integrand in the
cost functional (1.1) satisfies that
\begin{eqnarray}
F(x,u,t)\equiv C,\qquad\forall(x,u,t)\in \mathbb{R}^n\times U\times
[t_0,T],
\end{eqnarray}
for some constant $C\in \mathbb{R}$. Meanwhile, let
$f:\mathbb{R}^n\times U\times [t_0,T]\mapsto\mathbb{R}^n$ satisfies
the following assumption: Both $f$ and $f_x$ are continuous on
$\mathbb{R}^n\times U\times [t_0,T]$, and there exists a $L>0$ such
that
\begin{eqnarray}
\left\{\begin{array}{l}\|f(0,u,t)\|\le L,\\\|f_x(x,u,t)\|\le
L,\end{array}\right.\quad\forall (x,u,t)\in\mathbb{R}^n\times
U\times [t_0,T].
\end{eqnarray}

In these cases, any control $u(\cdot):[\tau, T]\mapsto U$ is an
optimal control under the initial data $(x^0,\tau)$, while the
Hamiliton-Jacobi-Bellman equation
\begin{eqnarray}
-v_{\tau}+\sup_{u\in U}[-\sum_{i=1}^nf_i(x,u,\tau)v_{x_i}-C]=0,
\end{eqnarray}
with the boundary condition
\begin{eqnarray}
v(x,T)=0,
\end{eqnarray}
admits a classical solution $V(x,\tau)=C(T-\tau)$, which is
obviously an extended solution of (3.3)-(3.4).

On the other hand, according to Pontryagin maximum principle
(\cite{PBGM}), any optimal control $u(\cdot)$ satisfies the
Pontryagin maximum condition.

Therefore, the Hamiliton-Jacobi-Bellman equation (3.3)-(3.4) has an
extended solution while more than one controls satisfies the
Pontryagin maximum condition.

\subsection{The second example}

Consider the special case of {\bf OCP} with $U=\mathbb{R}^m$ and the
quadratic cost functional, which is defined by
\begin{eqnarray}
F(x,u,t)=u^TSu,
\end{eqnarray}
where $S\in \mathbb{R}^{m\times m}$ is nonnegative semi-definite,
i.e., $u^TSu\ge 0$ for any $u\in \mathbb{R}^m$ and there exists
$u_0\ne 0$ such that $u_0^TSu_0=0$. Meanwhile, let
$f:\mathbb{R}^n\times U\times [t_0,T]\mapsto\mathbb{R}^n$ satisfies
the same assumption in the previous example.

In this case, the Hamiliton-Jacobi-Bellman equation
\begin{eqnarray}
-\frac{\partial v}{\partial \tau}+\sup_{u\in
U}[-\sum_{i=1}^nf_i(x,u,\tau)\frac{\partial v}{\partial
x_i}-u^TSu]=0,
\end{eqnarray}
with the boundary condition
\begin{eqnarray}
v(x,T)=0,
\end{eqnarray}
admits a classical solution $V(x,\tau)\equiv 0$, which is obviously
an extended solution of (3.6)-(3.7).

For any constant $k\in \mathbb{R}$, the control $u(t)=ku_0$ is an
optimal control under any initial data $(x_0,\tau)$. According to
Pontryagin maximum principle (\cite{PBGM}), any optimal control
$u(\cdot)$ satisfies the Pontryagin maximum condition.

In general, {\bf Hypothesis 2} may not hold for {\bf OCP}.

\section{On Hypotheses 3}

\subsection{Some Preparations}
Let $\phi:\mathbb{R}^n\times U\mapsto\mathbb{R}$ be a given
function.
\begin{eqnarray}
\Phi(x):=\sup_{u\in U}\{\phi(x,u)\},\qquad \forall x\in\mathbb{R}^n,
\end{eqnarray}
is an extended function, i.e., taking values in
$\mathbb{R}\cup\{+\infty\}$. For any $x\in\mathbb{R}^n$,
\begin{eqnarray}
M(x)=\arg\min_{u\in U}\phi(x,u):=\{u\in U|\hskip2mm
\phi(x,u)=\Phi(x)\},
\end{eqnarray}
which is possibly empty.
\begin{lemma} Assume that $\Phi(x)$ is finite for all
$x\in\mathbb{R}^n$, and there exists a modulus
$\omega:\mathbb{R}_+\times\mathbb{R}_+\mapsto\mathbb{R}_+$
satisfying $(2.4)$ such that,
\begin{eqnarray}
|\phi(x,u)-\phi(y,u)|\le \omega(R,\|x-y\|),\quad \forall x,y\in
B_R(\mathbb{R}^n),\quad u\in U,
\end{eqnarray}
where $B_R(\mathbb{R}^n)$ denotes the open ball in $\mathbb{R}^n$
with a radius $R>0$ centered at $0$.

Then it holds that
\begin{eqnarray}
\begin{array}{l}
\partial_-\Phi(x)\supseteq\partial_-\phi_x(x,u),\\
\partial_+\phi_x(x,u)\supseteq\partial_+\Phi(x),
\end{array}\qquad \text{ for all }u\in M(x).
\end{eqnarray}
\end{lemma}
{\bf Proof}\quad  Since $-\Phi(x)=\inf_{u\in U}[-\phi(x,u)]$,
applying Lemma 2.11 in \cite{BC} yields the
conclusions.\hfill$\Box$\vspace{2mm}

\begin{lemma} Assume that $U$ is compact, and $\phi:\mathbb{R}^n\times U\mapsto\mathbb{R}$
satisfies (4.3), and
\begin{description}
\item[(A)]  $\phi(\cdot,u)$ is differentiable at $x\in\mathbb{R}$ uniformly
in $u\in U$, i.e., there exists some modulus
$\omega_1:\mathbb{R}_+\mapsto\mathbb{R}_+$ with
\begin{eqnarray}
\lim_{R\rightarrow 0}\omega_1(R)=0,
\end{eqnarray}
such that
\begin{eqnarray}
|\phi(x+\Delta x,u)-\phi(x,u)-\frac{\partial\phi}{\partial
x}(x,u)\Delta x|\le\|\Delta x\|\omega_1(\|\Delta x\|),
\end{eqnarray}
for small $\Delta x$ and all $u\in U$;
\item[(B)] $\frac{\partial}{\partial
x}\phi(x,\cdot):U\mapsto\mathbb{R}^n$ is continuous;
\item[(C)] $\phi(x,\cdot):U\mapsto\mathbb{R}$ is lower semicontinuous.
\end{description}
Then $M(x)\ne\emptyset$, and
\begin{eqnarray}
\partial_-\Phi(x)=\overline{\rm co}Y(x),
\end{eqnarray}
where $Y(x):=\{\frac{\partial}{\partial x}\phi(x,u)|\hskip 2mm u\in
M(x)\}$;
\begin{eqnarray}
\partial_+\Phi(x)=\left\{\begin{array}{l}Y(x),\\\emptyset,\end{array}\right.\begin{array}{l}\text{if
}Y(x) \text{ is a singleton},\\\text{if }Y(x) \text{ is not a
singleton}.\end{array}
\end{eqnarray}
In particular, $\Phi$ is differentiable at $x$ if and only if $Y(x)$
is a singleton.

Moreover, $\Phi$ has the directional derivative in any direction
$v\in \mathbb{R}^n$, given by
\begin{eqnarray}
\frac{\partial\Phi}{\partial v}(x)=\max_{u\in
M(x)}\frac{\partial\phi}{\partial x}(x,u)v=\max_{p\in
\partial_-\Phi(x)}\langle p,v\rangle.
\end{eqnarray}
\end{lemma}
{\bf Proof}\quad  Since $-\Phi(x)=\inf_{u\in U}[-\phi(x,u)]$,
applying Proposition 2.13 in \cite{BC} yields the conclusions.
\hfill$\Box$\vspace{2mm}

\begin{definition} The function $H:\mathbb{R}^n\times\mathbb{R}^n\times
[t_0,T]\mapsto\mathbb{R}$ defined by
\begin{eqnarray}
H(x,p,t):=\sup_{u\in U}\mathscr{H}(x,p,u,t)\equiv\sup_{u\in
U}\sum_{i=1}^n\{p_if_i(x,u,t)-F(x,u,t)\},
\end{eqnarray}
is called the Hamiltonian related to {\bf OCP}, where the functions
$f=(f_1,\cdots, f_n)^T:\mathbb{R}^n\times U\times
[t_0,T]\mapsto\mathbb{R}^n$ and $F:\mathbb{R}^n\times U\times
[t_0,T]\mapsto\mathbb{R}$ are the data of {\bf OCP}.
\end{definition}

In this section, we will need  some technical assumptions on
$f:\mathbb{R}^n\times U\times [t_0,T]\mapsto \mathbb{R}^n$ and
$F:\mathbb{R}^n\times U\times [t_0,T]\mapsto \mathbb{R}$ as follows:

\begin{description}
\item[(H5)] Both $f$ and $F$ are continuous, and there exists an absolute modulus $\omega:\mathbb{R}_+\times\mathbb{R}_+\mapsto\mathbb{R}_+$
satisfying $(2.4)$ such that,
\begin{eqnarray}
|\phi(x,u)-\phi(y,u)|\le \omega(R,\|x-y\|),\quad \forall x,y\in
B_R(\mathbb{R}^n),\quad u\in U,
\end{eqnarray}
with $\phi:\mathbb{R}^n\times U\mapsto\mathbb{R}$ being
\begin{eqnarray}
\phi(x,u)=f_1(x,u,t),\hskip2mm\text{or}\hskip2mm
\cdots,\hskip2mm\text{or}\hskip2mm
f_n(x,u,t),\hskip2mm\text{or}\hskip2mm F(x,u,t),
\end{eqnarray}
for any given $t\in [t_0,T]$.
\item[(H6)] $\frac{\partial}{\partial
x}f:\mathbb{R}^n\times U\times [t_0,T]\mapsto \mathbb{R}^{n\times
n}$ and $\frac{\partial}{\partial x}F:\mathbb{R}^n\times U\times
[t_0,T]\mapsto \mathbb{R}^n$ are continuous, and there exists an
absolute modulus $\omega_1:\mathbb{R}_+\mapsto\mathbb{R}_+$
satisfying (4.5) such that $\phi:\mathbb{R}^n\times
U\mapsto\mathbb{R}$ in {\bf (H5)} satisfy (4.6), for any given $t\in
[t_0,T]$.
\end{description}

\begin{proposition} Assume that {\bf (H5)} holds, and both $f$ and $F$ are differentiable with respect to $x\in\mathbb{R}^n$.
Then it holds that
\begin{description}
\item[(I)]
$H(x,\cdot,t)$ is convex;
\item[(II)]
If
\begin{eqnarray}
A(x,p,t):=\{u\in U|\mathscr{H}(x,p,u,t)=H(x,p,t)\}\ne\emptyset,
\end{eqnarray}
and $H(x,\cdot,t)$ is differentiable at $p\in \mathbb{R}$, then
\begin{eqnarray}
f(x,u,t)\equiv c, \qquad\text{on}\qquad A(x,p,t);
\end{eqnarray}
If $A(x,p,t)\ne\emptyset$ and $H(\cdot,p,t)$ is differentiable at
$p\in \mathbb{R}$, then
\begin{eqnarray}
\frac{\partial}{\partial x}f(x,u,t)p-\frac{\partial}{\partial
x}F(x,u,t)\equiv c, \qquad\text{on}\qquad A(x,p,t).
\end{eqnarray}
\end{description}
\end{proposition}
{\bf Proof}\quad {\bf (I)} Let $\mathscr{H}(x,p,u,t)=\langle
f(x,u,t),p\rangle-f^0(x,u,t)$, we have
\begin{eqnarray*}
H(x,\lambda p_1+(1-\lambda)p_2,t)&=&\sup_{u\in
U}\mathscr{H}(x,\lambda
p_1+(1-\lambda)p_2,u,t)\\
&=&\sup_{u\in
U}[\lambda\mathscr{H}(x,p_1,u,t)+(1-\lambda)\mathscr{H}(x,p_2,u,t)]\\
&\le&\lambda H(x,p_1,t)+(1-\lambda)H(x,p_2,t),
\end{eqnarray*}
for any $\lambda\in [0,1]$.

{\bf (II)} Lemma 4.1 yields the conclusions.\hfill$\Box$\vspace{2mm}

\begin{corollary} Assume that $f$ and $F$ satisfy the assumptions in Proposition 4.1. If $H(\cdot,\cdot,t)$ is differentiable and
$A(x,p,t)\ne\emptyset$ at some
$(x,p)\in\mathbb{R}^n\times\mathbb{R}^n$. Then
\begin{eqnarray}
\left\{
\begin{array}{l}
\frac{\partial}{\partial p}H(x,p,t)=\{f(x,u,t)|\quad
\mathscr{H}(x,p,u,t)=H(x,p,t)\}\\\frac{\partial}{\partial
x}H(x,p,t)=\{\frac{\partial}{\partial
x}f(x,u,t)p-\frac{\partial}{\partial x}F(x,u,t)|\quad
\mathscr{H}(x,p,u,t)=H(x,p,t)\}.
\end{array}\right.
\end{eqnarray}
\end{corollary}
This corollary is just Proposition 3.2 in \cite{CF}.

\begin{proposition} Assume that {\bf (H5)-(H6)} holds, and the control domain $U$ is compact.
If
\begin{eqnarray}
A(x,p,t):=\{u\in U|\mathscr{H}(x,p,u,t)=H(x,p,t)\}\ne\emptyset,
\end{eqnarray}
and
\begin{eqnarray}
f(x,u,t)\equiv c, \qquad\text{on}\qquad A(x,p,t),
\end{eqnarray}
then $H(x,\cdot,t)$ is differentiable at $p\in \mathbb{R}$; If
$A(x,p,t)\ne\emptyset$ and
\begin{eqnarray}
\frac{\partial}{\partial x}f(x,u,t)p-\frac{\partial}{\partial
x}F(x,u,t)\equiv c, \qquad\text{on}\qquad A(x,p,t),
\end{eqnarray}
then $H(\cdot,p,t)$ is differentiable at $x\in \mathbb{R}$.
\end{proposition}
{\bf Proof}\quad Lemma 4.2 yields the
conclusions.\hfill$\Box$\vspace{2mm}

\subsection{Main Results}

It is well-known that, under some convex assumptions of the data $f$
and $F$, the necessary condition -- Pontryagin maximum principle is
also sufficient (see \cite{YZ, Ma} etc.). According to Theorem 2.5
in \cite{YZ}, we have
\begin{proposition} Assume that the control domain $U$ is convex and $\mathscr{H}(\cdot,p,\cdot,t):\mathbb{R}^n\times
U\mapsto\mathbb{R}$ is concave for all $(p,t)\in\mathbb{R}^n\times
[t_0,T]$. Then $u^*:[\tau,T]\mapsto U$ is an optimal control of {\bf
OCP} under the initial condition $x(\tau)=x^0$ with $(x^0,\tau)\in
\mathbb{R}^n\times [t_0,T)$, if and only if $u^*$ satisfies
Pontryagin maximum principle, i.e.,
\begin{eqnarray}
\mathscr{H}(x^*(t),p^*(t),u^*(t),t)=\max_{u\in
U}\mathscr{H}(x^*(t),p^*(t),u,t),\quad\text{a.e.}\quad t\in
[\tau,T],
\end{eqnarray}
where $(x^*(\cdot),p^*(\cdot))$ is the unique solution to the
control Hamiltonian system (1.5)-(1.6).

In particular, {\bf Hypotheses 3} is true.
\end{proposition}

\noindent{\bf Remark} It is obviously the fact that, {\bf Hypotheses
3} is true provided that there exists an optimal control under all
initial conditions. There exist many references on the existence of
optimal controls such as the famous Cesari' type conditions
\cite{Ce}, etc..
\vskip3mm

Now we will consider {\bf Hypotheses 3} directly through
Hamilton-Jacobi Theory based on some regularity condition of the
data $f$ and $F$ instead of the above convex conditions.

\begin{definition} Let $H:\mathbb{R}^n\times\mathbb{R}^n\times
[t_0,T]\mapsto\mathbb{R}$ be defined in Definition 4.1. The system
\begin{eqnarray}
\left\{\begin{array}{l}\dot{x}=\frac{\partial}{\partial
p}H(x,p,t),\\\dot{p}=-\frac{\partial}{\partial
x}H(x,p,t),\qquad\qquad p(T)=0,\end{array}\right.
\end{eqnarray}
is called the Hamiltonian system  related to {\bf OCP}.
\end{definition}
\begin{definition}
The Hamiltonian system (4.21) is called a complete system provided
that, for any $\xi\in \mathbb{R}^n$, (4.21) under the terminal
condition
\begin{eqnarray}
x(T)=\xi,
\end{eqnarray}
has a unique solution $(x(\cdot),p(\cdot))$ on $[t_0,T]$; Moreover,
for any sequence of solutions $(x^i(\cdot),p^i(\cdot))$ of (4.21)
with
$$
\lim_{i\rightarrow+\infty}(x^i(\tau_i),p^i(\tau_i),\tau_i)=(\xi,\eta,\tau)\in\mathbb{R}^n\times\mathbb{R}^n\times
[t_0,T],
$$
it holds that $(x^i(\cdot),p^i(\cdot))$ converge uniformly on
$[t_0,T]$ to the solution $(x(\cdot),p(\cdot))$ of (4.21) with
$(x(\tau),p(\tau))=(\xi,\eta)$.
\end{definition}
\begin{definition} The Hamiltonian system (4.21) related to {\bf OCP} is called to have a
shock at time $\tau\in [t_0,T)$ if there exist two solution
$(x^i,p^i):[t_0,T]\mapsto\mathbb{R}^n\times\mathbb{R}^n$ of (4.21)
with $i=1,2$, such that
\begin{eqnarray}
x^1(\tau)=x^2(\tau),\qquad p^1(\tau)\ne p^2(\tau).
\end{eqnarray}
\end{definition}

According to \cite{F} (see Ch.5, p.607-610) or \cite{CF}, it follows
that

\begin{lemma}
Assume that the Hamiltonian system (4.21) is a complete system and
also has no shock at all. If for every $\xi\in\mathbb{R}^n$, there
exists a control $u:[t_0,T]\mapsto U$ such that the solution
$(x(\cdot),p(\cdot))$ of (4.21) with the terminal condition
\begin{eqnarray}
x(T)=\xi,
\end{eqnarray}
satisfies
\begin{eqnarray}
\dot{x}(t)=f(x(t),u(t),t),\qquad\text{a.e}\quad t\in [t_0,T],
\end{eqnarray}
then $V:\mathbb{R}^n\times [t_0,T]\mapsto\mathbb{R}$ defined by
\begin{eqnarray}
V(x(\tau),\tau):=\int_{\tau}^TF(x(s),u(s),s)\,ds,
\end{eqnarray}
is the value function of {\bf OCP}.
\end{lemma}

\begin{theorem} Assume that {\bf (H5)-(H6)} holds and the Hamiltonian system (4.21) is a complete system.
If the control domain $U$ is a compact, then {\bf Hypotheses 3}
holds.
\end{theorem}
{\bf Proof}\quad  By {\bf Hypotheses 3}, there exists a unique
control $u^*: [\tau,T]\mapsto U$ satisfying Pontryagin maximum
principle under any initial data $(x^0,\tau)$. Thus
$\widetilde{V}:\mathbb{R}^n\times [t_0,T]\mapsto\mathbb{R}$ is
well-defined as follows:
\begin{eqnarray}
\widetilde{V}(x^0,\tau):=\int_{\tau}^TF(x^*(s),u^*(s),s)\,ds,
\end{eqnarray}
where $x^*(\cdot)$ is the solution of (1.2) with the initial
condition $x(\tau)=x^0$ and the control $u^*(\cdot)$. Meanwhile, it
follows from (4.16) and the assumptions that, the unique solution
$(x^*(\cdot),p^*(\cdot))$ of the control Hamiltonian system
(1.5)-(1.6) with $u(\cdot)=u^*(\cdot)$ is also the unique solution
of the Hamiltonian system (4.21) with the terminal condition
$$x(T)=x^*(T).$$

The assumptions and Lemma 4.3 tell us that we only need to prove the
Hamiltonian system (4.21) has no shock at all. Otherwise, if there
are two solutions $(x^i(\cdot),p^i(\cdot))$ of the Hamiltonian
system (4.21) with $i=1,2$, such that
\begin{eqnarray}
x^1(\tau)=x^2(\tau),\qquad p^1(\tau)\ne p^2(\tau),
\end{eqnarray}
for some $\tau\in [t_0,T)$, then it follows from (4.16) and the
famous Filippov's Lemma in \cite{Fi} (known as Measurable Selection
Theorem) that there exist two admissible controls $u^i:
[\tau,T]\mapsto U$ with $i=1,2$, ( i.e., both $u^1(\cdot)$ and
$u^2(\cdot)$ are Lebesque measurable ) such that both $u^1(\cdot)$
and $u^2(\cdot)$ satisfies the Pontryagin maximum condition together
with $(x^1(\cdot),p^1(\cdot))$ and $(x^2(\cdot),p^2(\cdot))$,
respectively, under the same initial data
$(x^0,\tau)=(x^1(\tau),\tau)$. This is contradictory to {\bf
Hypotheses 3}. The proof is completed. \hfill$\Box$\vspace{2mm}

\noindent{\bf Remark 4.2.}\quad In this paper, we take the
admissible control set at the initial time $\tau\in [t_0,T]$ as
$\cal{U}(\tau)$ defined in (2.12), i.e., all Lebesque measurable
functions on $[\tau,T]$. For other types of admissible control set
such as all piecewise continuous controls, the conclusions in this
paper are also valid.

\section{The conclusions}

In this paper, we study in detail three hypotheses on the optimal
control theory proposed by L. I. Rozonoer \cite{R}. {\bf Hypotheses
3} is only considered for the case with the smooth Hamiltonian. Now
we are considering the case with the non-smooth Hamiltonian.

\end{document}